\documentclass{ifacconf}

\usepackage[utf8]{inputenc}
\usepackage{natbib}        
\usepackage{amsmath,amsfonts,amssymb}
\usepackage{mathtools}
\usepackage{dsfont}
\usepackage[svgnames]{xcolor}
\usepackage{graphicx}
\usepackage[shortlabels]{enumitem}
\usepackage{float}
\usepackage{caption}
\usepackage{subcaption}
\usepackage{tikz}

\usepackage{microtype}
\graphicspath{{resources/}}

\newcommand{\rbl}{\left (}
\newcommand{\rbr}{\right )}

\newcommand{\al}{\left \langle}
\newcommand{\ar}{\right \rangle}

\newcommand{\nl}{\left\|}
\newcommand{\nr}{\right\|}
\newcommand{\cbl}{\left\lbrace }
\newcommand{\cbr}{\right\rbrace }

\newcommand{\Norm}[2][ ]{\nl #2 \nr_{#1}}
\newcommand{\SNorm}[1]{\Norm[\infty]{#1}}

\newcommand{\N}{\mathds{N}}

\newcommand{\R}{\mathds{R}}
\newcommand{\Rp}{\R_{\geq0}}

\newcommand{\eps}{\varepsilon}
\renewcommand{\phi}{\varphi}

\newcommand{\cC}{\mathcal{C}}

\newcommand{\cE}{\mathcal{E}}

\newcommand{\cG}{\mathcal{G}}

\newcommand{\cN}{\mathcal{N}}

\newcommand{\cT}{\mathcal{T}}

\newcommand{\oT}{\mathbf{T}}

\newcommand{\Impl}{\Longrightarrow}

\newcommand{\fa}{\forall \, }
\newcommand{\ex}{\exists \, }

\renewcommand{\d}{\mathrm{d}}

\newcommand{\setdef}[2]{\cbl#1\left|\vphantom{#1} #2\right.\cbr}

\newcommand{\dd}[2][ ]{\tfrac{\text{\normalfont d}#1}{\text{\normalfont d}#2}}

\makeatletter
\newcommand{\Itemlabel}[2]{#2\def\@currentlabel{#2}\label{#1}}
\makeatother

\makeatletter
\newcommand\setcurrentname[1]{\def\@currentlabelname{#1}}
\makeatother

\DeclareMathOperator*{\esssup}{ess\,sup}
\newcommand{\GL}{\text{GL}}

\DeclareMathOperator*{\graph}{graph}

\DeclareMathOperator*{\rf}{ref}
\DeclareMathOperator*{\loc}{loc}

\newcommand{\nocontentsline}[3]{}
\newcommand{\tocless}[2]{\bgroup\let\addcontentsline=\nocontentsline#1{#2}\egroup}
\newcommand{\toclesslab}[3]{\bgroup\let\addcontentsline=\nocontentsline#1{#2\label{#3}}\egroup} 

\newtheorem{definition}[thm]{Definition}

\begin{document}
\begin{frontmatter}

\title{Funnel control with input filter for nonlinear systems of relative degree two\thanksref{footnoteinfo}}
\thanks[footnoteinfo]{Funded by the Deutsche Forschungsgemeinschaft (DFG, German
Research Foundation) -- Project-ID 524064985.}
\author[Paderborn]{Dario Dennstädt}
\author[Paderborn]{Janina Schaa}
\author[Paderborn]{Thomas Berger}

\address[Paderborn]{Institut für Mathematik, Universit\"{a}t Paderborn, Warburger Stra\ss e~100, 33098~Paderborn, Germany
        (dario.dennstaedt@uni-paderborn.de, jschaa@mail.uni-paderborn.de, thomas.berger@math.upb.de). }
\begin{abstract}
We address the problem of output reference tracking for unknown nonlinear
multi-input, multi-output systems with relative degree two and bounded-input
bounded-state (BIBS) stable internal dynamics. We propose a novel model-free
adaptive controller that ensures the evolution of the tracking error within
prescribed performance funnel boundaries. By applying an output filter, the 
control objective is achieved without utilizing derivative information 
of system’s output. The controller is illustrated by a numerical example.
\end{abstract}


\end{frontmatter}

\section{Introduction}
Introduced in \cite{IlchRyan02a}, \emph{funnel control} is a well-established
adaptive high-gain control methodology for output-reference tracking of
nonlinear multi-input, multi-output systems. The controller operates under
minimal structural assumptions -- namely stable internal dynamics and a known
relative degree with a sign-definite high-frequency gain matrix  -- yet it
guarantees that the tracking error remains within prescribed performance bounds. This
framework offers robustness against disturbances and ensures transient
performance without relying on an explicit system model. As a result, funnel
control has been successfully applied to tracking problems in various domains,
such as DC-link power flow control~\citep{SenfPaug14}, control of industrial
servo-systems~\citep{Hackl17}, and temperature control of chemical reactor
models~\citep{IlchTren04}.

Recent research in funnel control has focused on overcoming practical
implementation challenges. Key areas of investigation include handling
measurement losses \citep{BergerMeasurement} and actuator
failures~\citep{Zhang2025},
multi-agent settings~\citep{Min22,Zhang2025multi}, 
input constraints \citep{Hu2022,BergerInput}, and
sampled-data \citep{Lanza24} or discrete-time systems \citep{Yun2023}, as well as
incorporating prior model knowledge to enhance performance \citep{BergerFMPC}.
For a comprehensive literature overview, we refer to the survey
by~\cite{BergIlch23}. Despite these advancements, a long-standing and persistent
difficulty for both funnel control and adaptive high-gain methods in general has
been the handling of systems with a relative degree exceeding one,
see~\cite{MORSE96}.

The prevalent approach to achieve prescribed output tracking for
higher-order systems relies on recursively defined error variables that
incorporate output derivatives, cf.~\cite{BergLe18a,BergIlch21}. However, obtaining
these derivatives in practice requires differentiating the output measurement,
which is an ill-posed problem, especially in the presence of noise \cite[Sec.
1.4.4]{hackl2012contributions}. This fundamental issue of derivative
availability has been addressed in several ways. 
Early results on tracking with prescribed accuracy, while only using output
information, for single-input, single-output linear time-invariant systems were
achieved by~\cite{Miller91}.
\cite{ILCHMANN2006396,IlchRyan07} employed a backstepping-based
filter to obtain a funnel control mechanism which uses output feedback only, 
but its design is complex and impractical due to high powers of a
typically large gain function. Alternatively, the \emph{funnel pre-compensator}
-- a high-gain observer structure -- was introduced by~\cite{BergerReis18} 
(with~\cite{lanza2022output} providing the feasibility proof for higher
relative degrees) to supply necessary auxiliary derivative
signals. While effective, this method also remains relatively complicated and
hard to implement, especially for systems with relative degree higher than two.
\cite{LANZA2024121} recently proposed a sampled-data funnel controller  
for control-affine systems of relative degree two that also avoids output derivatives. However, in order to derive a sufficiently fast sampling
rate and sufficiently large control signals, knowledge about the system to be
controlled is necessary; upper bounds of the system dynamics are assumed to be
known.

\emph{Prescribed performance control} -- a methodology related to funnel control
-- also ensures that a reference signal is tracked  within prescribed
boundaries, see e.g. \cite{BECHLIOULIS2008, BECHLIOULIS2014}. However, the
controller is designed for systems in a certain normal form and requires full
system state access. To overcome the need for full state information and rely
solely on output measurements, \cite{Dimanidis2020} employed a high-gain observer.
Similarly, \cite{Liu2021} and \cite{CHOWDHURY2019107} combined a high-gain observer with
the funnel control scheme. A drawback of the latter approaches is that the
observer parameters must be chosen ``sufficiently large'' a priori. A specific
lower bound, however, is not provided.

The present paper proposes a funnel controller that achieves prescribed output tracking
for nonlinear multi-input, multi-output systems of relative degree two. By
utilizing an output filter signal, the controller achieves its objective without
incorporating output derivatives. The proposed design is both simple and
has a lower complexity compared to previous methods. 
Furthermore, its reduced parameter set allows for more straightforward performance tuning.
Although this paper focuses on systems with relative degree two, the approach
lays a groundwork that we believe is generalizable to systems with
higher relative degree.

\subsubsection{Nomenclature:}
$\N$ and $\R$ denote natural and real numbers, respectively.
$\N_0:=\N\cup\{0\}$ and $\Rp:=[0,\infty)$.
$\Norm{x}:=\sqrt{\al x,x\ar}$~denotes the Euclidean norm of $x\in\R^n$.
$\cC^p(V,\R^n)$ is the linear space of $p$-times continuously  differentiable
functions $f:V\to\R^n$, where $V\subset\R^m$, and $p\in\N_0\cup \{\infty\}$.
$\cC(V,\R^n):=\cC^0(V,\R^n)$.
On an interval $I\subset\R$,  $L^\infty(I,\R^n)$ denotes the space of measurable and essentially bounded
functions $f: I\to\R^n$ with norm $\SNorm{f}:=\esssup_{t\in I}\Norm{f(t)}$,
$L^\infty_{\text{loc}}(I,\R^n)$ the set of measurable and locally essentially bounded functions.
Furthermore, $W^{k,\infty}(I,\R^n)$ is the Sobolev space of all $k$-times weakly differentiable functions
$f:I\to\R^n$ such that $f,\dots, f^{(k)}\in L^{\infty}(I,\R^n)$.

\section{Problem formulation}

\subsection{System class}

We consider second order nonlinear multi-input, multi-output systems of the form 
\begin{equation} \label{eq:Sys}
    \begin{aligned}
    \ddot{y}(t) &= R_1 y(t) + R_2 \dot{y}(t) + 
    f \big(\oT(y,\dot{y})(t)\big)+\Gamma u(t) \\
   y|_{[0,t_0]} &= y^0  \in \cC^{1}([0,t_0],\R^m)
    \end{aligned}
\end{equation}
with $t_0\ge 0$, initial trajectory $y^0$, 
control input~$u\in L^\infty_{\loc}([t_0,\infty), \R^m)$, 
and output ${y(t)\in\R^m}$ at time $t\geq 0$.
In the case of $t_0=0$, we identify $\cC^1([0,t_0],\R^m)$
with the vector space $\R^{2m}$. 
Then, the initial value condition in~\eqref{eq:Sys} is replaced by
${(y(t_0),\dot{y}(t_0))=(y^0_1,y^0_2)\in\R^{2m}}$.
Note that the control input $u$ and the system's output $y$ have the same dimension ${m\in\N}$.
The system consists of the \emph{unknown} matrices $R_1, R_2 \in \R^{m \times m}$, \emph{unknown} nonlinear function 
$f\in\cC(\R^q,\R^m)$, 
an \emph{unknown} positive definite matrix $\Gamma\in\GL_m(\R)$ 
and \emph{unknown} nonlinear operator ${\oT:\cC(\Rp,\R^m)\times\cC(\Rp,\R^m)\to L^\infty_{\loc}([t_0,\infty),\R^q)}$.
The operator~$\oT$ is causal, locally Lipschitz and satisfies a bounded-input bounded-output property. 
It is characterized in detail in the following definition.
\begin{definition} \label{Def:OperatorClass} 
For $n,q\in\N$, and $t_0\geq 0$, the set $\cT^{2n,q}_{t_0}$ denotes the class of operators $\oT:
\cC(\Rp,\R^n)\times\cC(\Rp,\R^n) \to L^\infty_{\loc} ([t_0,\infty), \R^{q})$
for which the following properties hold:
\begin{enumerate}[(i)]
    \item\emph{Causality}:  $\fa y_1,y_2\in\cC(\Rp,\R^n)^2$  $\fa t\geq t_0$:
    \[
        y_1\vert_{[0,t]} = y_2\vert_{[0,t]}
        \quad \Impl\quad
        \oT(y_1)\vert_{[t_0,t]}=\oT(y_2)\vert_{[t_0,t]}.
    \]
    \item\emph{Local Lipschitz}: 
    $\fa t \ge t_0 $ $\fa y \in \cC([0,t], \R^n)^2$ 
    ${\ex \Delta, \delta, c > 0}$ 
    $\fa y_1, y_2 \in \cC(\Rp, \R^n)^2$ with
    $y_1|_{[0,t]} = y_2|_{[0,t]} = y $ 
    and $\Norm{y_1(s) - y(t)} < \delta$,  $\Norm{y_2(s) - y(t)} < \delta $ for all $s \in [t,t+\Delta]$:
    \[
     \hspace*{-2mm}   \esssup_{\mathclap{s \in [t,t+\Delta]}}  \Norm{\oT(y_1)(s) \!-\! \oT(y_2)(s) }  
        \!\le\! c \ \sup_{\mathclap{s \in [t,t+\Delta]}}\ \Norm{y_1(s)\!-\! y_2(s)}\!.
    \] 
    \item\label{Item:PropBIBO}\emph{Bounded-input bounded-output (BIBO)}:
    $\fa c_0 > 0$ $\ex c_1>0$  $\fa y_1,y_2 \in \cC(\Rp, \R^n)$:
    \[
    \sup_{t \in \Rp} \Norm{y_1(t)} \le c_0 \ 
    \Impl \ \esssup_{t \in [t_0,\infty)} \Norm{\textbf{T}(y_1,y_2)(t)}  \le c_1.
    \]
\end{enumerate}
\end{definition}

\begin{rem}
    The BIBO property~\ref{Item:PropBIBO} of operator $\oT$
    allows us to conclude from the available information (the system output $y$)
    that the internal dynamics of the system (modeled by $\oT$) stay bounded.
    This was also assumed in previous works on pure output feedback funnel control, see, e.g.,
    \citep{BergerReis18,IlchRyan07,lanza2022output},
    but is a stronger property than the corresponding property 
   \begin{enumerate}
    \item[(iii')]\label{ItemPropBIBOAlt}
    $\fa c_0 > 0$ $\ex c_1>0$  $\fa z \in \cC(\Rp, \R^{2n})$:
    \[
    \sup_{t \in \Rp} \Norm{z(t)} \le c_0 \ 
    \Impl \ \esssup_{t \in [t_0,\infty)} \Norm{\textbf{T}(z)(t)}  \le c_1
    \]
    \end{enumerate}
    used in classical funnel control, where the output derivative information is
    assumed to be available, see, e.g., \cite{BergLe18a,BergIlch21}.
    The system \eqref{eq:Sys} explicitly contains linear terms as 
    they would be otherwise excluded by  property~\ref{Item:PropBIBO},
    contrary to property~(iii').
\end{rem}

\begin{definition} \label{Def:system-class}
    We say that the system~\eqref{eq:Sys} belongs to the system class
    $\cN^{m,2}$, written $(f,\Gamma, \oT,R_1,R_2) \in\cN^{m,2}$, if, for $t_0\geq0$ and
    some $q\in\N$, the following holds: $f\in\cC(\R^q,\R^m)$, $R_1, R_2 \in
    \R^{m \times m}$, $\Gamma\in\GL_m(\R)$  is a positive definite matrix, and
    ${\oT\in\cT^{2m,q}_{t_0}}$.
\end{definition}

\subsection{Control objective}\label{Sec:ContrObj}
\begin{figure}[ht]
\hspace{1cm}
\begin{center}
\begin{tikzpicture}[scale=0.43]
\tikzset{>=latex}
  \filldraw[color=gray!15] plot[smooth] coordinates {(0.15,4.7)(0.7,2.9)(4,0.6)(6,1.5)(9.5,0.6)(10,0.533)(10.01,0.531)(10.041,0.5) (10.041,-0.5)(10.01,-0.531)(10,-0.533)(9.5,-0.6)(6,-1.5)(4,-0.6)(0.7,-2.9)(0.15,-4.7)};
  \draw[thick] plot[smooth] coordinates {(0.15,4.7)(0.7,2.9)(4,0.6)(6,1.5)(9.5,0.6)(10,0.533)(10.01,0.531)(10.041,0.5)};
  \draw[thick] plot[smooth] coordinates {(10.041,-0.5)(10.01,-0.531)(10,-0.533)(9.5,-0.6)(6,-1.5)(4,-0.6)(0.7,-2.9)(0.15,-4.7)};
  \draw[thick,fill=lightgray] (0,0) ellipse (0.4 and 5);
  \draw[thick] (0,0) ellipse (0.1 and 0.533);
  \draw[thick,fill=gray!15] (10.041,0) ellipse (0.1 and 0.533);
  \draw[thick] plot[smooth] coordinates {(0,2)(2,-0.2)(4,0.4)(7,-0.25)(9,0.25)(10,0.15)};
  \draw[thick,->] (-2,0)--(12,0) node[right,above]{\normalsize$t$};
  \draw[thick,dashed](0,0.533)--(10,0.533);
  \draw[thick,dashed](0,-0.533)--(10,-0.533);
  \node [black] at (0,2) {\textbullet};
  \draw[->,thick](4,-3)node[right]{\normalsize$\inf_{s\geq0}\tfrac{1}{\phi(s)}$}--(2.5,-0.6);
  \draw[->,thick](3,3)node[right]{\normalsize$(0,e(0))$}--(0.07,2.07);
  \draw[->,thick](9,3)node[right]{\normalsize$1/\phi(t)$}--(7,1.4);
\end{tikzpicture}
\end{center}
\caption{Error evolution in a funnel $\mathcal F_{\varphi}$ with boundary $1/\varphi(t)$. The figure is based on~\cite[Fig.~1]{BergLe18a}, edited for present purpose.}
\label{Fig:funnel}
\end{figure}
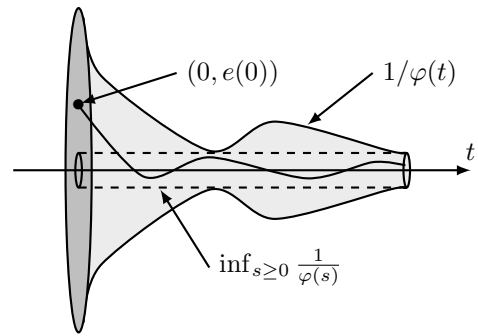

The objective is to design an output feedback control law,
which achieves that, for any reference signal ${y_{\rf}\in W^{2,\infty}(\Rp,\R^m)}$, 
the output tracking error ${e(t)=y(t)-y_{\rf}(t)}$ evolves within a prescribed performance funnel
\begin{equation*}
    \mathcal{F}_{\phi} := \setdef{(t,e)\in\R_{\ge 0} \times\R^m} {\phi(t)\|e\| < 1}.
\end{equation*}
This funnel is determined by the choice of the function~$\phi$ belonging to
\begin{align*}
    \cG:=\setdef
        {\varphi\in W^{1,\infty}(\Rp,\R)}
        {
          \begin{array}{l} \inf_{t\ge 0}  \phi(t) > 0
          \end{array}
        },
\end{align*}
see also Figure~\ref{Fig:funnel}.
Note that the tracking error $e$ evolving in $\mathcal{F}_{\phi}$ is not
constrained to converge to zero asymptotically. The choice of the performance
function~$\phi$ is guided by the application, which dictates the requisite
constraints on the tracking error.

\section{Controller design}\label{Sec:ContrDes}
We propose the following funnel controller to  achieve the control objective
described in Section~\ref{Sec:ContrObj} for systems of class $\cN^{m,2}$. 
\begin{equation}\label{eq:DefinitionFC}
\boxed{
\begin{aligned}
    e(t) &=y(t)-y_{\rf}(t),\\
    \dot{\xi}(t) &= -\xi(t)+u(t),& \xi(t_0)&=\xi^0,\\
    \theta(t)&=\xi(t)+\frac{e(t)}{1-\phi(t)^2\Norm{e(t)}^2},\\
    u(t)   &=\frac{-\theta(t)}{\hat{\theta}^2-\Norm{\theta(t)}^2},&\hat{\theta}&>0.
\end{aligned}
}    
\end{equation}
Note that the controller \emph{does neither} use the derivative~$\dot{y}$ 
of the system's output nor of the reference trajectory. Instead, we introduce a
filter variable $\xi$ that is designed to qualitatively reproduce $\dot{y}$. The
control law $u$ is then selected to ensure that $\xi$ tracks the virtual input signal
\[
    \xi^\star(t) = -\frac{e(t)}{1-\phi^2(t)\Norm{e(t)}},
\]
which is determined by the funnel controller for systems
of relative degree one, with the tracking error ${\theta(t)=\xi(t)-\xi^\star(t)}$ constrained to remain
in a funnel of constant radius $\hat{\theta} > 0$. In this way, the filter
variable acts as a surrogate for the unavailable derivative $\dot{y}$, enabling
the controller to achieve the desired control objective.
The controller has two parameters $\hat{\theta}$ and $\xi^0$, which determine the tracking performance.
However, note that no system knowledge is required to select the parameters, except for
 the initial system output $y^0(t_0)$, see Theorem~\ref{Thm:MainResult}.

A function $({x},\xi) : [0,\omega) \to \R^{2m}\times \R^m$, $\omega \in (t_0,\infty]$, is called solution of the closed-loop system~\eqref{eq:Sys},~\eqref{eq:DefinitionFC} in the sense of
\textit{Carath\'{e}odory}, if it satisfies the initial conditions $({x},\xi)|_{[0,t_0]} = (y^0,\dot y^0,\xi^0)$, and
 $({x},\xi)\vert_{[0,\omega)}$ is absolutely continuous and 
satisfies 
\begin{align*}
{\dot x}_1(t) &= {x}_{2}(t), \\
{\dot x}_2(t) &= R_1 x_1(t) + R_2 x_2(t) +  f \big(\oT(x)(t)\big)+\Gamma u(t)
\shortintertext{(which corresponds to~\eqref{eq:Sys} with $y=x_1$) and}
\dot{\xi}(t) &= -\xi(t)+u(t)
\end{align*}
 with~$u$ as in~\eqref{eq:DefinitionFC}  for almost all~$t\in[t_0,\omega)$.
We call a solution $x$ \emph{maximal}, if it does not have a right extension, which is also a solution.

We are now in the position to present the main result of the paper.
\begin{thm}\label{Thm:MainResult}
For $(f,\Gamma, \oT, R_1, R_2) \in\cN^{m,2}$, consider system~\eqref{eq:Sys}.
Let $\phi\in\cG, y_{\rf}\in W^{2,\infty}(\Rp,\R^m)$ and $\hat{\theta}>0$.
Further, let the initial trajectory $y^0  \in \cC^{1}([0,t_0],\R^m)$ be given with $\phi(t_0)^2\Norm{y^0(t_0)-y_{\rf}(t_0)}^2<1$ 
and choose $\xi^0\in\R^m$ such that 
\[
\Norm{\xi^0+\frac{y^0(t_0)-y_{\rf}(t_0)}{1-\phi(t_0)^2\Norm{y^0(t_0)-y_{\rf}(t_0)}^2}}<\hat{\theta}.
\]
Then, the application of the controller~\eqref{eq:DefinitionFC} 
to the system~\eqref{eq:Sys} yields an initial-value problem which has a solution, and 
every maximal solution $(x,\xi):[0,\omega)\to\R^{2m}\times\R^m$, $\omega\in(t_0,\infty]$, has the following properties:
\begin{enumerate}[(i)]
    \item\label{Item:GlobalSol} The solution is global, i.e. $\omega = \infty$.
    \item\label{Item:ControlBounded} The input $u:\Rp\to\R^m$ is bounded.
    \item\label{Item:FunnelProp} The tracking error~$e=x_1-y_{\rf}$ evolves uniformly within the performance funnel given by~$\phi$, i.e.
    \[
      \exists\, \eps\in(0,1)\  \forall\,t\geq t_0:\quad  \phi(t)\Norm{e(t)}<\eps.
    \]
\end{enumerate}
\end{thm}
\begin{pf}
\noindent
\emph{Step 1}:
    We show that there exists a maximal solution $(x,\xi):[0,\omega)\to\R^{2m}\times\R^m$, $\omega\in(t_0,\infty]$, of the closed-loop system~\eqref{eq:Sys},~\eqref{eq:DefinitionFC}.  Define the non-empty open set
    \begin{small}
    \[
        \cE\coloneqq\!\setdef{(t,z)\in\Rp\times\R^{3m}}
        {\!\!
        \begin{array}{l}
        \phi(t)\Norm{z_1-y_{\rf}(t)}<1,\\
             \Norm{z_3+\frac{z_1-y_{\rf}(t)}{1-\phi(t)^2\Norm{z_1-y_{\rf}(t)}^2}}<\hat{\theta} 
        \end{array}\!\!
        }\!.
    \] 
    \end{small}
    Using the notation $\tilde{\theta}(t,z_1,z_3)\coloneqq z_3+\frac{z_1-y_{\rf}(t)}{1-\phi(t)^2\Norm{z_1-y_{\rf}(t)}^2}$,
    define the function ${F:\cE\times\R^m\to\R^{m}}$ by 
    \begin{small}
    \[
        F(t,z_1,z_2,z_3,\eta) \coloneqq 
        \!
        \begin{bmatrix}
        z_2\\
        R_1z_1+R_2z_2+
        f(\eta)-\Gamma \frac{\tilde{\theta}(t,z_1,z_3)}{\hat{\theta}^2-\Norm{\tilde{\theta}(t,z_1,z_3)}^2}\\
        -z_3-\frac{\tilde{\theta}(t,z_1,z_3)}{\hat{\theta}^2-\Norm{\tilde{\theta}(t,z_1,z_3)}^2}
        \end{bmatrix}\!.
    \]
    \end{small}
    Then, the initial value problem~\eqref{eq:Sys},~\eqref{eq:DefinitionFC}
    takes the form 
    \[
        \begin{pmatrix}\dot{x}(t) \\ \dot \xi(t)\end{pmatrix}= F\left(t,\begin{pmatrix}{x}(t) \\  \xi(t)\end{pmatrix}, \oT (x)(t)\right),\quad \begin{pmatrix}{x}|_{[0,t_0]} \\  \xi|_{[0,t_0]}\end{pmatrix} = \begin{pmatrix} y^0 \\ \dot y^0 \\  \xi^0\end{pmatrix}.
    \]
    We have $(t_0,x(t_0),\xi^0)\in\cE$ by assumption.
    Application of a variant of~\cite[Thm.~B.1]{Ilchmann01102009} yields the existence of a maximal solution
    $(x,\xi):[t_0,\omega)\to\R^{2m}\times\R^m$, $\omega\in(t_0,\infty]$, with ${\graph((x,\xi)|_{[t_0,\omega)})\subset\cE}$.
    Moreover, the closure of $\graph\rbl (x,\xi)|_{[t_0,\omega)}\rbr$ is not a compact subset of~$\cE$.
    
    \noindent
    \emph{Step 2}: We define some constants for later use. Let $y(t):= x_1(t)$ for $t\in[t_0,\omega)$.
    By definition of the set $\cE$, we have $\phi(t)\Norm{y(t)}<1$ for all $t\in[t_0,\omega)$.
    On the entire interval $[t_0,\infty)$, 
    let $y^{e}$ be a continuous extension of $y$ with $\phi(t)\Norm{y^e(t)}<1$ 
    and $\tilde{y}^{e}$ be a continuous extension of $\dot{y}$ 
    (set $y^e\coloneqq y$ and $\tilde{y}^{e}\coloneqq \dot{y}$ in the case $\omega=\infty$).
    By the BIBO property of operator $\oT$, there exists a constant $c_1>0$ such that
    $\esssup_{t\in[t_0,\infty)}\Norm{\oT(y^e,\tilde{y}^{e})(t)}\leq c_1$.
    Therefore, ${(R_1\!-\!R_2\!-\!I) y^e(t)\!+\!f(\oT(y^e,\tilde{y}^{e})(t))}$ is bounded on the entire interval 
    $[t_0,\infty)$ due to the continuity of the involved functions. Consider the linear differential equation   
    \begin{equation}\label{Eq:AuxInternalDyn}
        \dot{\tilde{x}}(t)=-\tilde{x}(t)+(R_1\!-\!R_2\!-\!I)y^e(t) + f(\oT(y^e, \tilde{y}^e)(t))
    \end{equation}
    with initial value $\tilde{x}(t_0)=\dot{y}(t_0)-(R_2+I)y(t_0)-\Gamma\xi^0$. 
    Then,~\eqref{Eq:AuxInternalDyn} has a unique global solution $\tilde{x}:[t_0,\infty)\to\R^m$ that is bounded.
    Set $c_3\coloneqq \sup_{t\geq t_0}\Norm{\tilde{x}(t)}$
    and define the constants
    \begin{align*}
        \tilde{\lambda}&\coloneqq\Norm{R_2+I}\rbl\SNorm{y_{\rf}}+\SNorm{\tfrac{1}{\phi}}\rbr+c_3+\SNorm{\dot{y}_{\rf}}\\
        \shortintertext{and}
        \lambda&\coloneqq \SNorm{\tfrac{\dot{\phi}}{\phi}}+ \SNorm{\phi}\rbl\tilde{\lambda}+\Norm{\Gamma}\hat{\theta}\rbr,
    \end{align*}
    which are well defined by $\phi\in\cG$ and $y_{\rf}\in W^{2,\infty}(\Rp,\R^m)$. Due to the unboundedness of the function $x\mapsto\frac{x}{1-x}$ on $[0,1)$, it is possible to choose $\eps_1\in(0,1)$ such that
    ${\eps_1> \phi(t_0)\Norm{e(t_0)}}$ and 
    \[
        \frac{\eps_1^2}{1-\eps_1^2}> \frac{\lambda}{\lambda_{\min}(\Gamma)},
    \]
    where $\lambda_{\min}(\Gamma)>0$ denotes the smallest eigenvalue of the positive definite matrix $\Gamma$.
    
    \noindent
    \emph{Step 3}: We show $\phi(t)\|e(t)\|\leq\eps_1$ for all $t\in [t_0,\omega)$. 
    Seeking a contradiction, assume there exists $t^{\star}\in [t_0,\omega)$ with $\phi(t^\star)\Norm{e(t^\star)}>\eps_1$.
    As $\phi(t_0)\|e(t_0)\|<\eps_1$, there exists 
    \[
        t_{\star}:=\sup\setdef{t\in[t_0,t^\star)}{\phi(t)\|e(t)\|=\eps_1}<t^{\star},
    \]
    because $e$ is a continuous function on the entire interval $[t_0,t^\star]$.
    Using the shorthand notation 
    \[
        \xi^\star(t)\coloneqq-\frac{e(t)}{1 -\phi(t)^2\Norm{e(t)}^2},
    \]
    define the constants
    \[
        \bar{\xi}\coloneqq\sup_{s\in[t_0,t^\star]}\Norm{\xi(s)}\quad \text{and}\quad   \bar{\xi}^\star\coloneqq\sup_{s\in[t_0,t^\star]}\Norm{\dot{\xi}^\star(s)}
    \]
    and  choose $\eps_2\in(0,\hat{\theta})$ such that $\Norm{\theta(t_0)}<\eps_2$ and $\frac{\eps_2^2}{\hat{\theta}^2-\eps_2^2}>\hat{\theta}(\bar{\xi}+\bar{\xi}^\star)$.
    
    \noindent
    \emph{Step 3.1}: We show $\Norm{\theta(s)}\le\eps_2$ for all $s\in[t_0,t^\star]$ for $\theta$ as in~\eqref{eq:DefinitionFC}. 
    Seeking a contradiction, assume there exists $s^\star\in[t_0,t^\star]$ with $\Norm{\theta(s^\star)}>\eps_2$.
    Then, there exists 
    \[
        s_\star\coloneqq\sup\setdef{s\in[t_0,s^\star)}{\Norm{\theta(s)}=\eps_2}<s^\star.
    \]
    By construction of $s_\star$, we have $\Norm{\theta(s)}>\eps_2$ and
    $\frac{\Norm{\theta(s)}^2}{\hat{\theta}^2-\Norm{\theta(s)}^2}>\hat{\theta}(\bar{\xi}+\bar{\xi}^\star)$ for all $s\in(s_\star,s^\star]$.
    Omitting the dependency on $s$ and invoking~\eqref{eq:DefinitionFC}, we calculate
    \begin{align*}
         \tfrac{1}{2} \dd{s}\Norm{\theta}^2
        &=\al\theta,\dot{\theta}\ar=\al\theta,\dot{\xi}-\dot{\xi}^\star\ar\\
        &=\al\theta, -\dot{\xi}^\star\ar+\al\theta, -\xi\ar+\al\theta, u\ar\\
        &\le \hat{\theta}\rbl\bar{\xi}+\bar{\xi}^\star\rbr-\frac{\Norm{\theta}^2}{\hat{\theta}^2-\Norm{\theta}^2}<0
    \end{align*}
    almost everywhere on $[s_\star,s^\star]$, where we used $\|\theta(s)\|\le \hat \theta$, which holds by $(s,x(s),\xi(s))\in\mathcal{E}$. 
    Integration yields 
    \[
        \eps_2^2 \!<\! \| \theta(s^\star)\|^2\! =\!\!\int^{s^\star}_{s_\star}\!\!\!\!\dd{\tau}\Norm{\theta(\tau)}^2\d{\tau}
        +  \|\theta(s_\star)\|^2\le \|\theta(s_\star)\|^2 \!=\! \eps_2^2,
    \]
    a contradiction. Thus, $\Norm{\theta(s)}\le\eps_2$ for all $s\in[t_0,t^\star]$.

    \noindent
    \emph{Step 3.2}: We show that the existence of $t^\star$ leads to a contradiction.
    By definition of $t_\star$ and choice of $\eps_1$, we have 
    $1>\phi(t)\|e(t)\|\geq\eps_1$  and 
    \[
    \frac{\phi^2(t)\Norm{e(t)}^2}{1-\phi(t)^2\Norm{e(t)}^2}\ge \frac{\lambda}{\lambda_{\min}(\Gamma)}
    \]
    for all $t\in [t_\star,t^\star]$.
    Setting
    \begin{align}\label{eq:DefZeta}
        \zeta(t)\coloneqq \dot{y}(t)-(R_2+I)y(t)-\Gamma\xi(t),
    \end{align}
    we have  $\dot{y}(t)=(R_2+I) y(t)+ \zeta(t)+\Gamma\xi(t)$.
    Omitting the dependency on $t$, we calculate  
    \begin{align*}
        \dot{\zeta}&=\ddot{y}-(R_2+I)\dot{y}-\Gamma\dot{\xi}\\
        &=R_1y\!+\!R_2\dot{y}\!+\!f \big(\oT(y,\dot{y})\big)\!+\!\Gamma u\!-\!(R_2\!+\!I)\dot{y}\!-\!\Gamma(-\xi\!+\!u)\\
        &=R_1y+f \big(\oT(y,\dot{y})\big)-\dot{y}+\Gamma\xi\\
        &=-\zeta+(R_1-R_2-I)y+f \big(\oT(y,\dot{y})\big)
    \end{align*}
    on the interval $[t_0,\omega)$. Hence, $\zeta$ fulfills the differential equation~\eqref{Eq:AuxInternalDyn} 
    since $(y(t),\dot{y}(t))=(y^e(t),\tilde{y}^e(t))$ for all ${t\in[t_0,\omega)}$.
    Furthermore, the initial conditions coincide, $\zeta(t_0) = \tilde x(t_0)$,
    thus, by uniqueness of the solution of~\eqref{Eq:AuxInternalDyn}, $\zeta(t) = \tilde x(t)$ for all $t\in[t_0,\omega)$.
    Therefore, ${\Norm{\zeta(t)}\leq \Norm{\tilde{x}(t)}\leq c_3}$ for all $t\in[t_0,\omega)$.
    Omitting the dependency on $t$, we calculate 
    \begin{align*}
    & \tfrac{1}{2}\dd{t} \phi^2\Norm{ e}^2
    = \phi\dot{\phi}\Norm{e}^2+\phi^2\al e, \dot{e}\ar
    = \tfrac{\dot{\phi}}{\phi}\underbrace{\phi^2\Norm{e}^2}_{<1}+\phi^2\al e, \dot{e}\ar\\
    &\leq\SNorm{\tfrac{\dot{\phi}}{\phi}}+ \phi^2\al e,\dot{y}-\dot{y}_{\rf}\ar\\
    &\stackrel{\eqref{eq:DefZeta}}{=} \SNorm{\tfrac{\dot{\phi}}{\phi}}+ \phi^2\al e,(R_2+I)y+\zeta+\Gamma\xi-\dot{y}_{\rf}\ar\\
    &\leq\SNorm{\tfrac{\dot{\phi}}{\phi}}+ \SNorm{\phi}\Big(\Norm{R_2+I}\rbl\SNorm{y_{\rf}}+\SNorm{\tfrac{1}{\phi}}\rbr+c_3\\
    &\hspace{1.42cm}+\!\SNorm{\dot{y}_{\rf}}\Big)+\phi^2\al e, \Gamma\xi\ar\\
    &\le  \SNorm{\tfrac{\dot{\phi}}{\phi}}+ \SNorm{\phi}\tilde{\lambda}+\phi^2\al e, \Gamma(\theta+\xi^\star)\ar\\
    &\leq \SNorm{\tfrac{\dot{\phi}}{\phi}}+ \SNorm{\phi}\rbl\tilde{\lambda}+\Norm{\Gamma}\hat{\theta}\rbr+\phi^2\al e, -\frac{\Gamma e}{1-\phi^2\Norm{e}^2}\ar\\
    &\le \lambda-\lambda_{\min}(\Gamma)\frac{\phi^2\Norm{e}^2}{1-\phi^2\Norm{e}^2}\le 0 
    \end{align*}
   on $[t_\star,t^\star]$.
    Therefore,
    \[
        \eps_1^2 < \phi(t^\star)^2 \|e(t^\star)\|^2 
        \le \phi(t_\star)^2\|e(t_\star)\|^2 = \eps_1^2,
    \]
    a contradiction. Thus, we have  $\phi(t)\|e(t)\|\leq\eps_1$ for all $t\in [t_0,\omega)$.

    \noindent
    \emph{Step 4}: We show that $\xi$ is bounded on $[t_0,\omega)$.
    As a consequence of Step~3, we have  $\phi(t)\|e(t)\|\leq\eps_1$ for all ${t\in [t_0,\omega)}$. 
    Hence, there exists $\hat{\xi}^\star>0$ such that ${\|\xi^\star(t)\|\le \hat{\xi}^\star}$ for all $t\in [t_0,\omega)$.
    To show $\|\xi(t)\|\! \leq \!\max\{\|\xi(t_0)\|, \hat{\xi}^\star\}$ for all $t\in [t_0, \omega)$, we
    assume  that there exists $\hat{t}\in [t_0, \omega)$ such that $\|\xi(\hat{t})\|\! > \!\max\{\|\xi(t_0)\|, \hat{\xi}^\star\}$.
    Then, 
    \[ t_1\coloneqq\sup\setdef{s\in[t_0,\hat t)}{\|\xi(s)\| = \max\{\|\xi(t_0)\|, \hat{\xi}^\star\}} \] 
    is well-defined. Omitting the dependency on $t$,  we calculate
    \begin{align*}
        \tfrac{1}{2}\dd{t}\Norm{\xi}^2&=\al\xi,-\xi\ar+\al\xi,u\ar \\
        &= -\Norm{\xi}^2+\tfrac{1}{\hat{\theta}^2-\Norm{\theta}^2}\al\xi,-\theta\ar\\
        &= -\Norm{\xi}^2+\tfrac{1}{\hat{\theta}^2-\Norm{\theta}^2}\al\xi,-\xi+\xi^\star\ar\\
        &\leq -\Norm{\xi}^2+\tfrac{1}{\hat{\theta}^2-\Norm{\theta}^2}\rbl\Norm{\xi}\Norm{\xi^\star}-\Norm{\xi}^2\rbr\\
        &\leq -\Norm{\xi}^2+\tfrac{1}{\hat{\theta}^2-\Norm{\theta}^2}\Norm{\xi}\rbl\hat{\xi}^\star-\Norm{\xi}\rbr\leq 0
    \end{align*}
    on $[t_1,\hat{t}]$. Therefore,
    \begin{align*}
    \max\{\|\xi(t_0)\|, \hat{\xi}^\star\}^2\!&< \|\xi(\hat{t})\|^2 
    \le \|\xi(t_1)\|^2 \\
    &\leq \max\{\|\xi(t_0)\|, \hat{\xi}^\star\}^2,
    \end{align*}
    a contradiction. This shows that $\Norm{\xi(t)} \leq \max\{\|\xi(t_0)\|, \hat{\xi}^\star\}$ for all $t \in [t_0, \omega)$. Hence $\xi$ is bounded.
    
    \noindent
    \emph{Step 5}: We show that $\dot{y}$ is bounded on $[t_0,\omega)$. 
    As a consequence of Step~3, we have 
    $\phi(t)\|e(t)\|\leq\eps_1$ for all $t\in [t_0,\omega)$.
    Since~$y_{\rf}$ is bounded and $\inf_{t\geq t_0}\phi(t)>0$, the system output $y$ is bounded, too.
    Invoking~\eqref{eq:DefZeta}, we have
    \[
        \dot{y}(t)=(R_2+I) y(t)+\zeta(t)+\Gamma \xi(t)
    \]
    for all $t\in[t_0,\omega)$. Since $\zeta$ fulfills the differential equation~\eqref{Eq:AuxInternalDyn},
    we have ${\Norm{\zeta(t)}\leq c_3}$ for all $t\in[t_0,\omega)$, as shown in Step~3.2. According to 
    Step~4, the function $\xi$ is bounded. Therefore, $\dot{y}$ is bounded as well on the entire interval $[t_0,\omega)$.
    
    \noindent
    \emph{Step 6}: We show that $u$ is bounded.
    According to the Steps~3--5, the functions $y$, $\dot{y}$, and $\xi$ are bounded functions and $\varphi \|e\|$ is uniformly bounded away from~1.
    Therefore, $\xi^\star$ is a bounded function as well.  
    Since $y_{\rf}\in W^{2,\infty}(\Rp,\R^m)$, we find that $\dot \xi^\star$ is bounded 
    and hence we may define the constants
    \[
        \bar{\xi}\coloneqq\sup_{s\in[t_0,\infty)}\Norm{\xi(s)}\quad \text{and}\quad   \bar{\xi}^\star\coloneqq\sup_{s\in[t_0,\infty)}\Norm{\dot{\xi}^\star(s)}.
    \]
    By choosing $\eps_3\in(0,\hat{\theta})$ such that $\Norm{\theta(t_0)}<\eps_3$ and 
    ${\frac{\eps_3^2}{\hat{\theta}^2-\eps_3^2}>\hat{\theta}(\bar{\xi}+\bar{\xi}^\star)}$,
    we can adapt Step~3.1 to show ${\Norm{\theta(t)}\le \eps_3}$ for all $t\in[t_0,\infty)$.
    As a consequence, $u$ is bounded. 
    
    \noindent
    \emph{Step 6}: We show $\omega=\infty$.
    According to Step~5, $\dot{y}$ is a bounded function. Thus,
    there exists $\eps_4>0$ with $\SNorm{\dot{y}}\leq \eps_4$.
    According to Steps~3--6, $\graph\rbl (x,\xi)|_{[t_0,\omega)}\rbr$ is
    contained in the set
    \begin{small}
    \[
    \setdef{(t,z)\in\Rp\times\R^{3m}}
        {\!\!
        \begin{array}{l}
        \phi(t)\Norm{z_1-y_{\rf}(t)}\leq\eps_1,\\
        \Norm{z_2}\leq\eps_4,\\
             \Norm{z_3+\frac{z_1-y_{\rf}(t)}{1-\phi(t)^2\Norm{z_1-y_{\rf}(t)}^2}}\!\!\leq\eps_3\!
        \end{array}\!\!
        }\!\subset \cE 
    \]
    \end{small}
    Since the closure of $\graph\rbl x|_{[t_0,\omega)}\rbr$
    is not a compact subset of $\cE$ according to the observation from Step~1, 
    this implies $\omega=\infty$ and thereby shows assertion~\ref{Item:GlobalSol}.
    Further, $\phi(t)\|e(t)\|\leq \eps_1<1$ for all $t\in[t_0,\infty)$  shows assertion~\ref{Item:FunnelProp}.
    Moreover, $u$ is bounded on the whole interval $[t_0,\infty)$ according to Step~6
    verifying assertion~\ref{Item:ControlBounded}.
    This completes the proof.\hfill $\Box$
\end{pf}

\section{Simulations}

\begin{figure}
\includegraphics[width = \columnwidth]{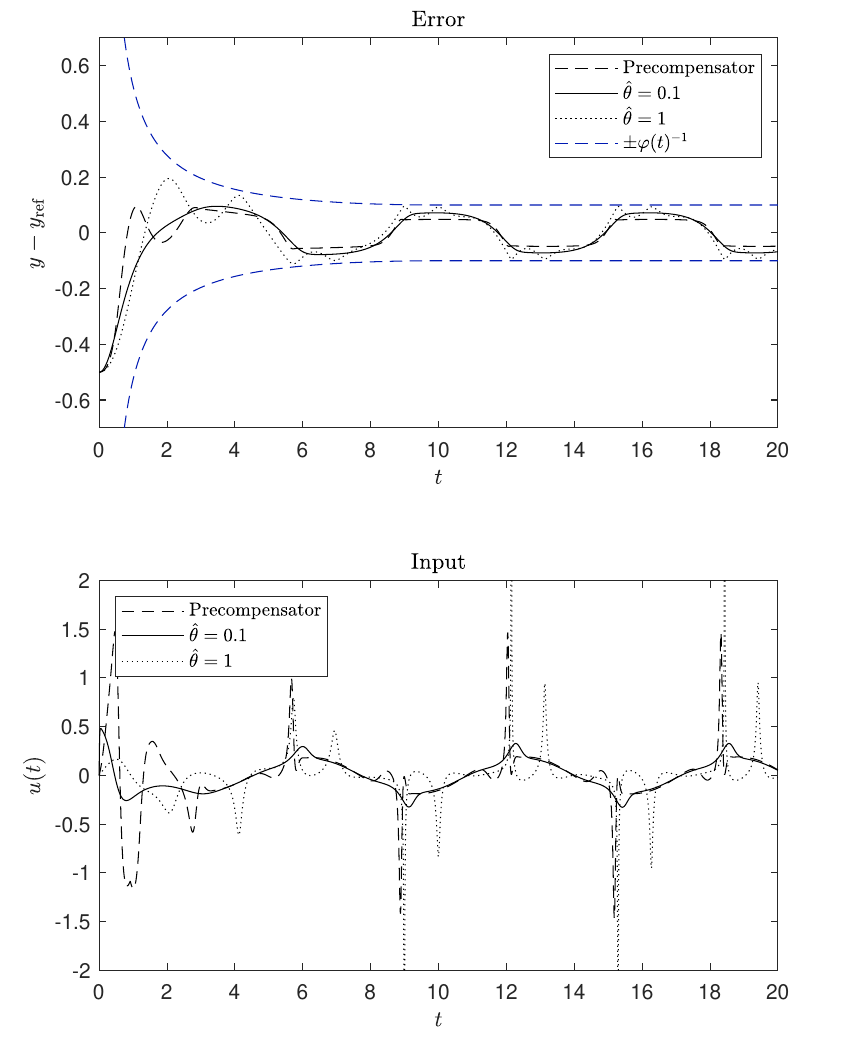}
\caption{Simulation of the behavior of the error $y(t) - y_{\text{ref}}(t)$ and the input signal generated by the controller proposed by~\cite{8062785} and the controller~\eqref{eq:DefinitionFC}, for the system~\eqref{eq:sys-sim}, using the parameters $a = 2, b = 2$. 
}
\label{Fig:Simulation}
\end{figure}

To illustrate the previous results, we simulate the controller~\eqref{eq:DefinitionFC} for different choices of the parameter $\hat{\theta}$ and compare it to the performance of the controller proposed by~\cite{8062785}. This controller introduces two auxiliary variables $z_1, z_2$, yielding the following structure: 
\begin{align*}
\dot{z}_1(t) &= z_2(t)+(q_1+p_1k_2(t)) (y(t)-z_1(t)),\\
\dot{z}_2(t) &= (q_2+p_2k_2(t))(y(t)-z_1(t))+\tilde{\Gamma}u(t),\\
e_0(t) &= z_1(t) - y_{\text{ref}}(t),\\
e_1(t) &= \dot{e}_0(t) + k_0(t)e_0(t), \text{ where} \\
k_0(t) &= \frac{1}{1-\varphi_0(t)^2 \Vert e_0(t)\Vert^2},\\
k_1(t) &= \frac{1}{1-\varphi_1(t)^2 \Vert e_1(t)\Vert^2},\\
k_2(t) &= \frac{1}{1-\varphi_2(t)^2 \Vert y(t)-z_1(t) \Vert^2},\\
u(t) &= -k_1(t)e_1(t),
\end{align*}
for a suitable choice of constants $p_1, p_2, q_1, q_2 \in \R$, ${\tilde{\Gamma}
\in \R^{m \times m}}$ and performance functions $\varphi_0, \varphi_1,
\varphi_2$. Evidently, this controller involves much more design parameters,
which need to be adjusted appropriately, than the controller proposed
in~\eqref{eq:DefinitionFC}. The latter also exhibits a lower complexity, as only
one differential equation is involved. In order to compare the performance of
both controllers, we consider the system also examined by \cite{8062785} given by 
\begin{equation}\label{eq:sys-sim}
    \ddot{y}(t) + a \sin y(t) = bu(t),
\end{equation}
with the system parameters $a = 2, b = 2$. The simulations are
performed in MATLAB (solver: \texttt{ode45}, rel.tol.: $10^{-8}$, abs.tol.:
$10^{-6}$) over the time interval $[0,20]$. The results for tracking the reference trajectory $y_{\text{ref}}(t) = \frac{\cos(t)}{2}$ with the initial state $(y(0),\dot{y}(0)) = (0,0)$ for the choices $\xi^0 = 0.5$ for the initial value of the filter variable and $z_1(0) = z_2(0) = 0.5$ for the auxiliary variables for the controller from~\cite{8062785} are depicted in
Figure~\ref{Fig:Simulation}. The performance functions ${\varphi_1(t) = \frac{1}{e^{-t}+1}}$ and ${\varphi_0(t) = \varphi_2(t) = 20 \begin{cases} 1-(0.1t-1)^2, & t \leq 10, \\ 1, & t > 10, \end{cases}}$ and design parameters $\tilde{\Gamma} = 2$, $p_1 = 1$, $p_2 = \frac{5}{7}, q_1 = 1, q_2 = 5$ are obtained by sequential tuning of the parameters, while the funnel function ${\varphi(t) = (\varphi_0(t)^{-1} + \varphi_2(t)^{-1})^{-1}}$ is chosen such that it ensures the same transient performance of the error $y(t) - y_{\text{ref}}(t)$. 

We observe that the performance of the proposed controller~\eqref{eq:DefinitionFC} depends on the choice of parameter $\hat{\theta}$. 
For sufficiently small values of $\hat{\theta}$ (smaller than $0.1$), further reductions do not noticeably alter the system response, which yields a straightforward tuning process. Since it is much more challenging to properly tune the controller proposed in~\cite{8062785} due to the larger number of design parameters, finding a well-performing configuration is considerably more difficult. The input signal does not show any peaks, and has a smaller bandwidth compared to the controller proposed by~\cite{8062785}. Also, the output shows improved tracking behavior. However, for larger values of $\hat{\theta}$, the input signal shows larger bandwidth and substantial peaks, while the output tracking is inferior. This highlights that an appropriate selection of $\hat{\theta}$ is crucial for satisfactory controller performance.\\


\section{Conclusion}

In this work, we developed a funnel controller for systems of relative degree
two that achieves trajectory tracking with prescribed performance and without requiring derivatives of
the output signal by introducing a filter variable. The controller structure is
simpler than in previous approaches, and the number of design parameters is much smaller, allowing for a more direct influence of the controller performance by appropriate selection, which is also exhibited in simulations.
Future work will focus on extending this approach to systems of
higher relative degree. 
\small
\bibliography{references}
\end{document}